\documentclass[12pt]{article}

\newcommand{\qed}{{\hfill\rule{4pt}{7pt}}}

\usepackage{amsmath, epsfig, cite, lineno}
\usepackage{amssymb}
\usepackage{amsfonts}
\usepackage{latexsym}

\newtheorem{thm}{Theorem}[section]

\newtheorem{cor}[thm]{Corollary}

\newtheorem{lem}[thm]{Lemma}
\newtheorem{conj}[thm]{Conjecture}

\def\pf{\noindent {\it Proof.} }

\numberwithin{equation}{section}

\begin{document}


\begin{center}
{\Large\bf Proof of Sun's conjecture on the divisibility of certain binomial sums}
\end{center}

\vskip 2mm \centerline{Victor J. W. Guo}
\begin{center}
{\footnotesize Department of Mathematics, East China Normal
University,\\ Shanghai 200062,
 People's Republic of China\\
{\tt jwguo@math.ecnu.edu.cn,\quad
http://math.ecnu.edu.cn/\textasciitilde{jwguo} }}
\end{center}

\vskip 0.7cm
{\small \noindent{\bf Abstract.} In this paper, we prove the following result conjectured by Z.-W. Sun: 
$$
(2n-1){3n\choose n}\left|
\sum_{k=0}^{n}{6k\choose 3k}{3k\choose k}{6(n-k)\choose 3(n-k)}{3(n-k)\choose n-k}\right.
$$
by showing that the left-hand side divides each summand on the right-hand side. 
}



\section{Introduction}
In \cite{Sun}, Z.-W. Sun proved some new series for $1/\pi$ as well as related congruences on sums
of binomial coefficients, such as
\begin{align*}
\sum_{n=0}^\infty\frac{n}{864^n}\sum_{k=0}^{n}{6k\choose 3k}{3k\choose k}
{6(n-k)\choose 3(n-k)}{3(n-k)\choose n-k}
=\frac{1}{\pi},
\end{align*}
and, for any prime $p>3$,
\begin{align*}
\sum_{n=0}^{p-1}\frac{n}{864^n}\sum_{k=0}^{n}{6k\choose 3k}{3k\choose k}
{6(n-k)\choose 3(n-k)}{3(n-k)\choose n-k}
\equiv 0\pmod{p^2}.
\end{align*}
Sun \cite{Sun} also proposed many interesting related conjectures, one of which is 
\begin{conj}\label{conj:sun}
{\rm\cite[Conjecture 4.2]{Sun}}For $n=0,1,2,\ldots$ define
\begin{align*}
s_n:=\frac{1}{(2n-1){3n\choose n}}
\sum_{k=0}^{n}{6k\choose 3k}{3k\choose k}{6(n-k)\choose 3(n-k)}{3(n-k)\choose n-k}.
\end{align*}
Then $s_n\in\mathbb{Z}$ for all $n$. Also,
\begin{align}
\lim_{n\to\infty}\sqrt[n]{s_n}=64.  \label{eq:64}
\end{align}
\end{conj}
Sun himself has proved that $s_n\equiv 0\pmod 8$ for $n\geq 1$ and
$s_{p-1}\equiv \lfloor(p-1)/6 \rfloor\pmod p$ for any prime $p$.
In this paper, we shall prove that Conjecture {\rm\ref{conj:sun}} is true by
establishing the following two theorems.
\begin{thm}\label{thm:1}
For $0\leq k\leq n$, there holds
\begin{align*}
\frac{1}{(2n-1){3n\choose n}}
{6k\choose 3k}{3k\choose k}{6(n-k)\choose 3(n-k)}{3(n-k)\choose n-k}\in\mathbb{Z}.
\end{align*}
\end{thm}
Note that, in \cite{Sun2,Sun3}, Sun proved many similar results on the divisibility
of binomial coefficients.
\begin{thm}\label{thm:2}
For $n\geq 2$ and $0\leq k< n/2$, there holds
\begin{align*}
&\hskip -2mm
{6k\choose 3k}{3k\choose k}{6(n-k)\choose 3(n-k)}{3(n-k)\choose n-k} \\
&\geq {6k+6\choose 3k+3}{3k+3\choose k+1}{6(n-k-1)\choose 3(n-k-1)}{3(n-k-1)\choose n-k-1},
\end{align*}
and so
\begin{align}
\frac{2}{2n-1}{6n\choose 3n}\leq s_n \leq \frac{n+1}{2n-1}{6n\choose 3n}. \label{eq:2n-1}
\end{align}
\end{thm}
It is easy to see that \eqref{eq:64} follows from \eqref{eq:2n-1} and Stirling's formula
$$
n!\sim \sqrt{2\pi n} \left(\frac{n}{e}\right)^n.
$$

\section{Proof of Theorem \ref{thm:1}}
We need the following two lemmas.
\begin{lem} Let $m,n\geq 1$ and $0\leq k\leq n$. Then 
\begin{align*}
{mn\choose n}\left|{2mk\choose mk}{mk\choose k}{2m(n-k)\choose m(n-k)}{m(n-k)\choose n-k}.\right.
\end{align*}
\end{lem}
\pf Observe that
\begin{align}
&\hskip -2mm
{2mk\choose mk}{mk\choose k}{2m(n-k)\choose m(n-k)}{m(n-k)\choose n-k}\left/{mn\choose n}\right. \nonumber \\
&=\frac{(2mk)!(2mn-2mk)!}{(mk)!(mn-mk)!(mn)!}{(m-1)n\choose (m-1)k}{n\choose k}.  \label{eq:mkmk}
\end{align} 
The proof then follows from the fact that numbers of the form
$$
\frac{(2a)!(2b)!}{a!b!(a+b)!},
$$
called the {\it super Catalan numbers}, are integers (see \cite{Gessel,Warnaar}).
\qed
\begin{lem}\label{lem:2}
Let $0\leq k\leq n$ be integers. Then
\begin{align*}
(2n-1)\left|\frac{(6k)!(6n-6k)!(2n)!}{(3k)!(3n-3k)!(3n)!(2k)!(2n-2k)!}\right..
\end{align*}
Or equivalently,
\begin{align*}
\frac{(6k)!(6n-6k)!(2n)!(2n-2)!}{(3k)!(3n-3k)!(3n)!(2k)!(2n-2k)!(2n-1)!}\in\mathbb{Z}.
\end{align*}
\end{lem}
\begin{cor}Let $n\geq 1$. Then $(2n-1)$ divides ${6n\choose 3n}$.
\end{cor}
The order in which a prime $p$ enters $n!$ can be written as
\begin{align}
{\rm ord}_pn!=\sum_{i=1}^\infty\left\lfloor\frac{n}{p^i}\right\rfloor,  \label{eq:ord}
\end{align}
where $\lfloor x\rfloor$ denotes the greatest integer less than or equal to $x$.
In order to prove Lemma \ref{lem:2}, we first establish the following result.
\begin{lem}Let $m\geq 2$ and $0\leq k\leq n$ be integers.
Then
\begin{align}
&\hskip -2mm \left\lfloor\frac{6k}{m}\right\rfloor+\left\lfloor\frac{6n-6k}{m}\right\rfloor
+\left\lfloor\frac{2n}{m}\right\rfloor+\left\lfloor\frac{2n-2}{m}\right\rfloor \nonumber\\
&\geq \left\lfloor\frac{3k}{m}\right\rfloor+\left\lfloor\frac{3n-3k}{m}\right\rfloor
+\left\lfloor\frac{3n}{m}\right\rfloor+\left\lfloor\frac{2k}{m}\right\rfloor
+\left\lfloor\frac{2n-2k}{m}\right\rfloor+\left\lfloor\frac{2n-1}{m}\right\rfloor, \label{eq:long}
\end{align}
unless $k\not\equiv 1\pmod 3$, $m=3$, and $n\equiv 2\pmod 3$.
\end{lem}
\pf For any real numbers $x$ and $y$, it is well known that
\begin{align*}
\lfloor 2x\rfloor +\lfloor 2y\rfloor &\geq \lfloor x\rfloor +\lfloor y\rfloor+\lfloor x+y\rfloor, \\
\lfloor x+y\rfloor &\geq \lfloor x\rfloor +\lfloor y\rfloor.
\end{align*}
It follows that
\begin{align*}
\left\lfloor\frac{6k}{m}\right\rfloor+\left\lfloor\frac{6n-6k}{m}\right\rfloor
+\left\lfloor\frac{2n}{m}\right\rfloor 
\geq\left\lfloor\frac{3k}{m}\right\rfloor+\left\lfloor\frac{3n-3k}{m}\right\rfloor
+\left\lfloor\frac{3n}{m}\right\rfloor+\left\lfloor\frac{2k}{m}\right\rfloor
+\left\lfloor\frac{2n-2k}{m}\right\rfloor.
\end{align*}
Suppose that \eqref{eq:long} does not hold. Then we must have
\begin{align}
\left\lfloor\frac{6k}{m}\right\rfloor+\left\lfloor\frac{6n-6k}{m}\right\rfloor
&=\left\lfloor\frac{3k}{m}\right\rfloor+\left\lfloor\frac{3n-3k}{m}\right\rfloor
+\left\lfloor\frac{3n}{m}\right\rfloor, \label{eq:3.1}\\
\left\lfloor\frac{2n}{m}\right\rfloor 
&=\left\lfloor\frac{2k}{m}\right\rfloor+\left\lfloor\frac{2n-2k}{m}\right\rfloor, \label{eq:3.2}\\ 
\left\lfloor\frac{2n-2}{m}\right\rfloor
&<\left\lfloor\frac{2n-1}{m}\right\rfloor,  \label{eq:3.3}
\end{align}
and so, by \eqref{eq:3.3}, $m\mid 2n-1$, then by \eqref{eq:3.2}, $m\mid 2k$ or $m\mid 2n-2k$. 
If $m\mid 2k$, then $m\mid k$ (since $m\mid 2n-1$ means that $m$ is odd), 
the identity  \eqref{eq:3.1} implies that
\begin{align}
\left\lfloor\frac{6n}{m}\right\rfloor
=2\left\lfloor\frac{3n}{m}\right\rfloor. \label{eq:6nm}
\end{align}
Since $m\mid 2n-1$, the identity \eqref{eq:6nm} implies that
\begin{align}
\frac{2n-1}{m}+\left\lfloor\frac{3}{m}\right\rfloor
=2\left\lfloor\frac{n+1}{m}\right\rfloor. \label{eq:2n-1m}
\end{align}
If $m\geq 4$, then the left-hand side of \eqref{eq:2n-1m} equals $(2n-1)/m$,
while the right-hand side of \eqref{eq:2n-1m} belongs to 
$$\left\{\frac{2n+2}{m},\,\frac{2n}{m},\,\frac{2n-2}{m},\ldots\right\},$$
a contradiction, and so $m\leq 3$. Since $m\geq 2$ is odd, we must have $m=3$. 
Hence, $n\equiv 2\pmod 3$ and $k\equiv 0\pmod 3$. Similarly, if $m\mid 2n-2k$, then we deduce that
$m=3$, $n\equiv 2\pmod 3$ and $k\equiv 2\pmod 3$. This proves the lemma.
\qed

\medskip
\noindent{\it Proof of Lemma \ref{lem:2}.} 
For any prime $p\neq 3$, by \eqref{eq:ord} and \eqref{eq:long}, we have 
\begin{align}
&{\rm ord}_p(6k)!+{\rm ord}_p(6n-6k)!+{\rm ord}_p(2n)!+{\rm ord}_p(2n-2)! \nonumber \\
&\geq {\rm ord}_p(3k)!+{\rm ord}_p(3n-3k)!+{\rm ord}_p(3n)!
+{\rm ord}_p(2k)!+{\rm ord}_p(2n-2k)!+{\rm ord}_p(2n-1)!. \label{eq:ineq}
\end{align}
For $p=3$, since ${\rm ord}_3(3j)!=j+{\rm ord}_3 j!$, the inequality \eqref{eq:ineq} 
reduces to
\begin{align}
{\rm ord}_3(2n)!+{\rm ord}_3(2n-2)! 
\geq {\rm ord}_3 k!+{\rm ord}_3(n-k)!+{\rm ord}_3 n!+{\rm ord}_3(2n-1)!. \label{eq:p3}
\end{align}
Noticing that
$$ 
\frac{(2n)!(2n-2)!}{k!(n-k)! n!(2n-1)!}
=\frac{1}{2n-1}{2n\choose n}{n\choose k}
=\left(4{2n-2\choose n-1}-{2n\choose n}\right){n\choose k}\in\mathbb{Z},
$$
the inequality \eqref{eq:p3} holds. Namely, the inequality \eqref{eq:ineq}
is true for $p=3$. This completes the proof.  \qed

\medskip
\noindent{\it Proof of Theorem \ref{thm:1}.} By  \eqref{eq:mkmk}, one sees that
\begin{align*}
&\hskip -2mm \frac{1}{(2n-1){3n\choose n}}
{6k\choose 3k}{3k\choose k}{6(n-k)\choose 3(n-k)}{3(n-k)\choose n-k}  \\
&=\frac{1}{2n-1}\frac{(6k)!(6n-6k)!}{(3k)!(3n-3k)!(3n)!}{2n\choose 2k}{n\choose k},
\end{align*}
which is an integer divisible by ${n\choose k}$ in view of Lemma \ref{lem:2}.  \qed

\section{Proof of Theorem \ref{thm:2}}
Let 
\begin{align*}
A_{n,k}={6k\choose 3k}{3k\choose k}{6(n-k)\choose 3(n-k)}{3(n-k)\choose n-k}.
\end{align*}
Then, for $n\geq 2$ and $0\leq k<n/2$, we have
\begin{align*}
\frac{A_{n,k}}{A_{n,k+1}}-1
=\frac{(36nk+31n-36k^2-36k-5)(n-2k-1)}{(6k+5)(6k+1)(n-k)^2}\geq 0,
\end{align*}
i.e., $A_{n,k}\geq A_{n,k+1}$. Since $A_{n,k}=A_{n,n-k}$, for $n\geq 2$, we have
$$
2A_{n,0}=2{6n\choose 3n}{3n\choose n}\leq \sum_{k=0}^n A_{n,k} \leq (n+1)A_{n,0}.
=(n+1){6n\choose 3n}{3n\choose n}
$$
In other words, the inequality \eqref{eq:2n-1} holds. 

\renewcommand{\baselinestretch}{1}


\begin{thebibliography}{99}
\small \setlength{\itemsep}{-.8mm}
\bibitem{Gessel}I. Gessel, Super ballot numbers. J. Symb. Comput. 14 (1992), 179--194.

\bibitem{Sun}Z.-W. Sun, Some new series for $1/\pi$ and related congruences, preprint, arXiv:1104.3856.

\bibitem{Sun2}Z.-W. Sun, Products and sums divisible by central binomial coefficients,
¡¡Electron. J. Combin. 20(1) (2013), \#P9. 

\bibitem{Sun3}Z.-W. Sun, On divisibility of binomial coefficients, J. Austral. Math. Soc., in press. 

\bibitem{Warnaar}S.O. Warnnar, A $q$-rious positivity, Aequat. Math. 81 (2011), 177--183.

\end{thebibliography}
\end{document}